\documentclass[a4paper,11pt]{article}

\usepackage{latexsym,amssymb,amsmath,cite,epsfig} 
\usepackage{a4wide}   
\usepackage[notref,notcite]{}
\usepackage{framed}
\usepackage{float}
\usepackage{hyperref}
\usepackage{alltt}

\long\def\del#1\enddel{}

\def\fBM#1\fEM{\hbox{\fns$\begin{matrix}#1\end{matrix}$}}
\def\BE {\begin{equation}}      \def\EE {\end{equation}}	
\let\fns=\footnotesize		
\def\IP{{\mathbb P}}		 
\let\ni=\noindent

\begin{document}          \baselineskip=16pt

\thispagestyle{empty}
\begin{flushright}
TUW-11-12
\end{flushright}
\vspace{1cm}
\begin{center}
{\LARGE A new offspring of PALP}
\end{center}
\vspace{8mm}
\begin{center}
{\large Andreas~P.~Braun\footnote{e-mail: abraun AT hep.itp.tuwien.ac.at} and Nils-Ole~Walliser\footnote{e-mail: walliser AT hep.itp.tuwien.ac.at}}
\end{center}
\vspace{3mm}
\begin{center}
{\it Institut f\"ur Theoretische Physik, Vienna University of Technology,\\
Wiedner Hauptstrasse 8-10, 1040 Vienna, Austria
}
\end{center}
\begin{center}
\end{center}
\vspace{15mm}
\begin{abstract}
\noindent We describe the C program \textit{mori.x}. It is part of PALP, a package for analyzing lattice polytopes. Its main purpose is the construction 
and analysis of three--dimensional smooth Calabi--Yau hypersurfaces in toric varieties. The ambient toric varieties are given in terms of fans over the 
facets of reflexive lattice polytopes. The program performs crepant star triangulations of reflexive polytopes and determines the Mori cones of the resulting 
toric varieties. Furthermore, it computes the intersection rings and characteristic classes of hypersurfaces.
\end{abstract}
\vspace{2cm}
\begin{flushright} 
{\it In memory of Maximilian Kreuzer}
\end{flushright} 
\vspace{1.2cm}
\newpage
\setcounter{tocdepth}{2}
\tableofcontents
\setcounter{footnote}{0}


\section{Introduction}
PALP \cite{Kreuzer:2002uu}, a package for analyzing lattice polytopes, has the ability to construct, manipulate and analyze
lattice polytopes. The present program \textit{mori.x}, which is also included in the new version of PALP (starting from release 2.0) 
available at \cite{PALPwww}, adds further functionalities concerning Calabi--Yau (CY) three--fold hypersurfaces.
The smoothness of the CY hypersurface is achieved by triangulating the reflexive polytope, i.e. appropriately resolving 
the ambient space. This functionality is optimized for the CY hypersurface case. Generic CY hypersurfaces avoid point--like 
singularities of the ambient space as well as divisors that correspond to interior points of facets.
Consequently, our algorithm performs crepant star triangulations only up to such interior points.\footnote{Complete triangulations of 
arbitrary polytopes can be performed with the program TOPCOM \cite{topcom}, which is also included in the open--source mathematics software 
system SAGE \cite{sage}. SAGE also contains various tools for handling toric varieties.}  Using these triangulations, the SR--ideals and 
the intersection rings are computed. The latter are determined with the help of SINGULAR \cite{DGPS}, a computer algebra system 
for polynomial computations. Furthermore, we have implemented the Oda--Park algorithm \cite{Odapark, Berglund:1995gd} to find the Mori cone 
of the ambient space.

The program can also analyze arbitrary three--dimensional hypersurfaces embedded in the ambient toric varieties.
It is capable of computing the intersection ring, characteristic classes and Hodge numbers.  
As we do not triangulate the polyhedron completely, non--CY hypersurfaces can also contain singularities.
Hence, the intersection numbers of non--CY hypersurfaces determined by \textit{mori.x} are not integers in general.

The triangulations, which are at the heart of this program, are subject to two types of limitations. First, \textit{mori.x} performs 
crepant star triangulation of the polytope by subdividing its facets into simplices. This procedure is, however, only implemented 
for those non--simplicial cones which give rise to at most three--dimensional secondary fans. Second, the limitations of PALP apply,
as \textit{mori.x} is a part of this package. PALP works with fixed precision; maximal dimensions, maximal number of polytope points, 
vertices and facets need to be fixed before compiling. By editing the file \textit{Global.h}, the user can change the default values if needed. 
See the PALP documentation \cite{Kreuzer:2002uu} for further details.

This note is meant as a manual for \textit{mori.x}, an introduction to the mathematical background is beyond its scope.
We refer the reader to \cite{fulton,generaltoric} for introductions to toric geometry and the construction of CY hypersurfaces
using polytopes. The structure of the paper is as follows. In section~2 we explain the I/O formats, present the options of \textit{mori.x}
and describe the basic functions of the main program. This section should be used as a quick reference
for the various functionalities. Aspects which deserve further explanations are presented in section~3. There, we also discuss the 
resolution of singular compact toric spaces via (crepant star) triangulations. Section~4 describes the structure 
of the program. The principal routines are listed and commented on in the header file \textit{Mori.h}.
\\\\
\ni\textbf{Acknowledgments:}
We dedicate this work to Maximilian Kreuzer who passed away before it was completed. He contributed the core routines to the source code of \textit{mori.x}. 
The present publication is an occasion to savor his work. We hope that we have been able to realize the original spirit in which he started this project.  

We would like to thank Harald Skarke for valuable advice in the completion phase of the program, and Christoph Mayrhofer for his comments and contributions to 
the source code at an early stage. Furthermore, we would like to thank Johanna Knapp, Volker Braun, Maximilian Attems, Ching--Ming Chen and Michal Michal$\check{\text{c}}$ik  
for discussion and comments.

The work of N.-O.~Walliser was supported by the FWF under grant P21239 and I192.
The work of A.~P.~Braun was supported by the FWF under grant P21239 and P22000. 

\section{I/O formats and options of \textit{mori.x}}

In this section we explain the I/O formats of \textit{mori.x} and give an overview of its options. This program is part of the latest version 
of the package PALP that is available at the web page \cite{PALPwww} under the GNU license terms. There, the reader can find a detailed installation 
description. In short, type ``make all" in the command line to compile the entire suite of programs. Otherwise,  type ``make mori" to compile 
\textit{mori.x} only. Consult \cite{Kreuzer:2002uu} for the documentation of the other functions of the package: \textit{poly.x}, \textit{cws.x}, 
\textit{class.x} and \textit{nef.x}.

Let us start with the help screen. It provides essential information about all the functionalities of the
program (letters in boldface denote the command line input):

{\fns
\begin{alltt}
\textbf{\$ mori.x -h}
This is ``mori.x'':  performing crepant star triangulations of a polytope P* in the N-lattice
                     computing the Mori cone of the corresponding toric ambient spaces
                     determining the intersection rings of embedded (CY) hypersurfaces
Usage:   mori.x [-<Option-string>] [in-file [out-file]]
Options (concatenate any number of them into <Option-string>):
    -h      print this information 
    -f      use as filter
    -g      general output: triangulation and Stanley-Reisner ideal
    -I      incidence information of the facets (ignoring interior points of facets)
    -m      Mori generators of the ambient space
    -P      IP-simplices among points of P* (ignoring interior points of facets)
    -K      points of P* in Kreuzer polynomial form
    -b      Hodge numbers and Euler number
            or arithmetic genera and Euler number if combined with -H
    -i      intersection ring
    -c      Chern classes of the (CY) hypersurface
    -t      triple intersection numbers
    -d      topological information on toric divisors and del Pezzo conditions
    -a      all of the above except h, f, I and K
    -D      lattice polytope points of P* as input (default CWS)
    -H      arbitrary hypersurface class `H = c1*D1 + c2*D2 + ...' as input (default CY)
            Input: coefficients `c1 c2 ...'
Input: 1) standard input format: degrees and weights `d1 w11 w12 ... d2 w21 w22 ...'
       2) alternative (use -D): `d np' or `np d' (d=Dimension, np=#[points])
                                and (after newline) np*d coordinates
Output:   as specified by options

\end{alltt}
}

As the program works with reflexive polytopes only, the input of the program must specify such a polytope. There are two
input formats available: one either gives a combined weight system (CWS) or directly provides a matrix of lattice points 
whose convex hull is the polytope. The input is checked for reflexivity, see \cite{Kreuzer:2000qv, Kreuzer:2002uu} for 
details on how a CWS is converted into a polytope. 

The construction of a CY hypersurface uses a dual pair of reflexive lattice polytopes. One polytope determines the CY hypersurface,
whereas its dual specifies the ambient toric variety. In this work, we are mainly interested in the $N$--lattice polytope which gives 
rise to the fan of the ambient toric variety. We follow the conventions of the literature and PALP in referring to this polytope 
as the dual polytope $P^*$.

Let us consider a $\IP^1$ fibered over a $\IP^3$ as an example:

{\fns\begin{alltt}
\textbf{\$ mori.x -P}
Degrees and weights  `d1 w11 w12 ... d2 w21 w22 ...':
5 1 1 1 1 1 0  2 0 0 0 0 1 1
4 7  points of P* and IP-simplices
    1    0    0    0   -1    0    0
    0    1    0    0   -1    0    0
    0    0    1    0   -1    0    0
    0    0    0    1    1   -1    0
------------------------------   #IP-simp=2
    1    1    1    0    1    1   5=d  codim=0
    0    0    0    1    0    1   2=d  codim=3
\end{alltt}}
\ni Note that the ordering of the CWS input is not obeyed by the output of lattice points. Once the order is displayed, however, it
is fixed and determines the labeling of toric divisors in any further output.  

An alternative way to provide the input is to type lattice polytope points directly. In this case, one has to use the parameter ``-D". 
Let us reconsider the example above:

{\fns\begin{alltt}
\textbf{\$ mori.x -DP}
`#lines #colums' (= `PolyDim #Points' or `#Points PolyDim'):
\textbf{4 7}
Type the 28 coordinates as dim=4 lines with #pts=7 columns:
\textbf{1 0 0 -1  0 0 0
	0 1 0 -1  0 0 0
	0 0 1 -1  0 0 0
	0 0 0  1 -1 1 0}
4 7  points of P* and IP-simplices 
    1    0    0   -1    0    0    0
    0    1    0   -1    0    0    0
    0    0    1   -1    0    0    0
    0    0    0    1   -1    1    0
------------------------------   #IP-simp=2  
    1    1    1    1    1    0   5=d  codim=0
    0    0    0    0    1    1   2=d  codim=3
\end{alltt}}
\ni The order of the lattice points displayed is the same as the order set in the input. 
This type of input should be preferred if one wants to control the order of the polytope points.
Note that the IP--simplices among points of an $M$--lattice polytope $P$ can be obtained from the lattice points of $P^*$ with the option ``-P'' of \textit{poly.x}.
Since we focus on the $N$--lattice polytope $P^*$, this option is suppressed in \textit{mori.x} for the sake of simplicity.

The rest of this section contains a detailed description of the options listed in the help screen. If no flag is specified, the 
program starts with the parameter ``-g". By default, the program considers a CY hypersurface embedded in the
ambient toric variety. The option ``-H'' has to be used in order to consider non--CY hypersurfaces.

\begin{description}	\def\Item[#1]{\item[\hbox to 22 pt{#1\hss}]}

\Item[-h] This option prints the help screen.

\Item[-f] This parameter suppresses the prompt of the command line. This is useful if one wants to build pipelines.

\Item[-g] First, the triangulation data of the facets is displayed. The number of triangulated simplices is followed by the incidence 
structure of the simplices. The incidence information for each simplex is encoded in terms of a binary number: there is a digit for each
polytope point; a $1$ denotes that the point belongs to the simplex. Second, the Stanley--Reisner ideal is displayed: the number of 
elements of the ideal is followed by its elements. Each element is denoted by a binary number as above. 

\Item[-I] The incidence structure of the facets of the polytope $P^*$ is displayed. Interior points of the facets are neglected.

\Item[-m] The Mori cone generators of the ambient space are displayed in the form of a matrix. Each row corresponds to a generator. 
The entries of each row are the intersections of the generator with the toric divisor classes. The Oda--Park algorithm is used
to compute the generators.  

\Item[-P] Those IP--simplices whose vertices are lattice points of the polytope, but not inner points of the facets, are displayed. 

\Item[-K]  The Kreuzer polynomial is displayed, see section \ref{kreuzerpoly} for details. The number of points in the interior of facets is shown
as \verb+intpts+. The multiplicities of the toric divisors are displayed as \verb+multd+ if they are greater than one. Furthermore, the Picard 
number of the CY hypersurface is computed and printed as \verb+Pic+. 

\Item[-b] We distinguish two cases. For the Calabi--Yau hypersurface (default) case, the Hodge numbers $h^{1,1}$, $h^{2,1}$ and the Euler 
characteristic are displayed. If the input is an arbitrary hypersurface (see option ``-H"), the zeroth and first arithmetic genera are 
displayed instead of the Hodge numbers. As a non--CY hypersurface can be singular, care is needed for an interpretation of the
results; see section~\ref{sect31} for more details.

\Item[-i] This option displays the intersection polynomial restricted to a CY hypersurface. The polynomial 
is displayed in terms of an integral basis of the toric divisors. The coefficients of the monomials are the triple intersection 
numbers in this basis. This option can also be used together with ``-H" to perform this task for non--CY hypersurfaces. 

\Item[-c] The Chern classes of the hypersurface (CY or non--CY) are displayed in terms of an integral basis of the toric divisors. 

\Item[-t] The triple intersection numbers of the toric divisors are displayed. When the intersection numbers are equal to zero, they are not 
shown.

\Item[-d] This option displays topological data of the toric divisors restricted to the (CY or non--CY) hypersurface. The Euler characteristics of the toric 
divisor classes and their arithmetic genera are shown. The toric divisor classes are tested against necessary conditions for 
del Pezzo surfaces. The following data is listed: the del Pezzo candidates preceded by their number (their type is given in 
brackets) and those del Pezzo candidates that do not intersect other del Pezzos.

\Item[-a] This is a shortcut for ``-gmPbictd".

\Item[-D] This tells the program to expect a matrix of lattice points of the polytope $P^*$ as the input. This is useful if one wants to set and control the order of the toric divisors.

\Item[-H] Using this option, one can specify a (non--CY) hypersurface. The user determines the hypersurface divisor class $H = \sum_i c_i D_i$ in terms of the toric divisor classes $D_i$
by typing its coefficients $c_i$. The hypersurface can then be analyzed by combining ``-H'' with other
options, as described above. Just using ``-H'', the program runs ``-Hb''. See section~\ref{sect31} for an example.

\end{description}

\noindent
Note that the options ``-b, -i, -c, -t, -d, -a, -H'' need SINGULAR to be installed.


\section{Supplementary details}

We designed this program with an eye for applications to string theory model building, in particular the analysis of CY three--folds. 
In this context, our motivation was to implement the construction of \textit{smooth} CY hypersurfaces embedded in toric varieties, using the
construction of Batyrev \cite{Batyrev:1994hm}. Here, the starting point is a dual pair of reflexive polytopes which determines both the ambient
toric variety and a CY hypersurface.

\subsection{Triangulations and point--like singularities}\label{sect31}

Batyrev \cite{Batyrev:1994hm} has shown that any four--dimensional reflexive polytope gives rise to a smooth CY hypersurface after triangulation. 
The ambient toric space is not necessarily smooth, as some of the cones might have a volume greater than one even after triangulation. This leads to 
point--like singularities, which, however, do not meet a generic CY hypersurface. The reflexivity of the polytopes sensibly 
simplifies the triangulation procedure. The induced simplices of lower dimension do not contribute any further singularities.

\textit{mori.x} performs crepant star triangulations of a polyhedron ignoring points in the interior of facets. 
Even though this introduces further point--like singularities into the ambient toric variety, these are also avoided by 
generic CY hypersurfaces.

Polytopes can be triangulated by subdividing the secondary fans of its non--simplicial facets \cite{Billera, GKZ}. 
This triangulation algorithm is implemented in \textit{mori.x} for up to three--dimensional secondary fans. The program exits with 
a warning message if the subdivision is not properly completed.

Consider the following CWS:

{\fns \begin{alltt}
\textbf{\$  echo -e '8 4 1 1 1 1 0  6 3 1 0 1 0 1' | mori.x -fPI}
4 8  points of P* and IP-simplices
    3    1    0    0    0   -1    1    0
   -1    0    1    0    0    0    0    0
    3    0    0    0    1   -1    1    0
   -4    0    0    1    0    1   -1    0
------------------------------   #IP-simp=2  
    1    1    1    0    1    4   8=d  codim=0
    1    0    1    1    0    3   6=d  codim=1
Incidence: 101011 001111 111110 110101 011101 111001 100111
\end{alltt}}
\noindent This system describes a four--dimensional lattice polytope with eight points, six of which are vertices. We label the 
column vectors with $v_1,\dots,v_8$. Point $v_7$ lies in the interior of the third facet. The incidence data show the intersections 
of the six polytope points with the seven facets. The third facet contains the five points $v_1,\dots,v_5$, hence it is not simplicial and 
we have to triangulate it.

{\fns \begin{alltt}
\textbf{\$ echo -e '8 4 1 1 1 1 0  6 3 1 0 1 0 1' | mori.x -fg}
9 Triangulation
101011 001111 110101 011101 111001 100111 011110 110110 111010
2 SR-ideal
101100 010011
8 Triangulation
101011 001111 110101 011101 111001 100111 101110 111100
2 SR-ideal
010010 101101
\end{alltt}}
\noindent The program performs the two possible crepant triangulations of the facet $\langle 12345 \rangle$. 
The first result yields the three simplices $\langle \widehat{1 3 4} 2 5\rangle$, whereas the second gives the two simplices $\langle 1 \widehat{25}34\rangle$
(in this notation the hat indicates that one of the points is dropped). Nevertheless, the two resolutions give the same CY intersection polynomial:

{\fns\begin{alltt}
\textbf{\$ echo -e '8 4 1 1 1 1 0  6 3 1 0 1 0 1' | mori.x -fi}
SINGULAR -> divisor classes (integral basis J1 ... J2):
d1=J1, d2=-3*J1+J2, d3=J1, d4=4*J1-J2, d5=-3*J1+J2, d6=J2
SINGULAR -> intersection polynomial:
2*J1^3+108*J2^3+8*J1^2*J2+30*J2^2*J1 
SINGULAR -> divisor classes (integral basis J1 ... J2):
d1=J1, d2=-3*J1+J2, d3=J1, d4=4*J1-J2, d5=-3*J1+J2, d6=J2
SINGULAR -> intersection polynomial:
2*J1^3+108*J2^3+8*J1^2*J2+30*J2^2*J1 
\end{alltt}}

\noindent \verb+d1+,$\,\dots$, \verb+d6+ denote the toric divisors corresponding to the lattice points $v_1,\dots,v_6$. There are two independent divisor 
classes. Indeed, \textit{mori.x} expresses the intersection polynomial in terms of the integral basis $J_1=D_1=D_3$ and $J_2=D_6$.

Let us take a closer look at non--CY hypersurfaces. The reader is warned: for these cases smoothness is not guaranteed anymore, so that
the intersection numbers can become fractional. Some choices of the hypersurface equation may intersect point--like singularities not resolved by the triangulation. 
Consider e.g. the hypersurface divisor class $H =  D_1 + D_6$. Remember that the order in which \textit{mori.x} expects the coefficients of the hypersurface 
divisor class is fixed by the polytope matrix and not by the CWS input. Hence, the correct input for $H$ is the string \verb+1 0 0 0 0 1+.  

{\fns\begin{alltt}
\textbf{\$ mori.x -H}
Degrees and weights  `d1 w11 w12 ... d2 w21 w22 ...'
8 4 1 1 1 1 0  6 3 1 0 1 0 1
Type the 6 (integer) entries for the hypersurface class:
1 0 0 0 0 1
Hypersurface degrees: ( 5  4 )
Hypersurface class: 1*d1 1*d6 
SINGULAR  -> Arithmetic genera and Euler number of H:
chi_0: 35/32 , chi_1: 143/32  [ -27/4 ]
SINGULAR  -> Arithmetic genera and Euler number of H:
chi_0: 29/27 , chi_1: 128/27  [ -22/3 ]
\end{alltt}}
\noindent To calculate these quantities, the program determines the characteristic classes of the divisors using adjunction. 
It then performs the appropriate integration with the help of the triple intersection numbers. The fractional results of the arithmetic 
genera and the Euler number in our example indicate that the intersection polynomial has fractional entries. This happens because the 
singularity of the ambient toric variety descends to the hypersurface $H$. In the first triangulation, the simplex $\langle 1235\rangle$ still 
has volume four. In the second triangulation, there are two simplices with volume three: $\langle 1 2 3 4\rangle$ and $\langle 1 3 4 5\rangle$. 
Indeed, any hypersurface in the divisor class $H =  D_1 + D_6$ is forced to pass through the corresponding singularities.

\subsection{The Kreuzer polynomial}\label{kreuzerpoly}

The Kreuzer polynomial encodes lattice polytope points in a compact form. The number of variables equals the
dimension of the polytope. Each lattice point gives rise to a Laurent monomial in which the exponents of the 
variables are the coordinates. Vertices and non--vertices are distinguished by coefficients 
``+" and ``-" respectively. Points in the interior of facets are ignored. Consider the example presented above.

{\fns \begin{alltt}
\textbf{\$ echo -e '8 4 1 1 1 1 0  6 3 1 0 1 0 1' | ./mori.x -fK}
KreuzerPoly=t_1^3t_3^3/(t_2t_4^4)+t_1+t_2+t_4+t_3+t_4/(t_1t_3); intpts=1;  Pic=2
\end{alltt}}
\noindent Note that negative coordinates are always displayed by putting the variables in the denominator.

\subsection{The Mori cone of the ambient space}

For toric varieties, \textit{mori.x} uses the algorithm of Oda and Park \cite{Odapark, Berglund:1995gd} to compute the Mori cone. 
The generators of the Mori cone are given in terms of their intersections with the toric divisors. For singular
toric varieties, the Picard group of Cartier divisors is a non--trivial subgroup of the Chow group, which contains the Weil
divisors. Hence one can consider the K\"ahler cone, which is dual to the Mori cone, as a cone in the vector space
spanned by the elements of either the Picard or the Chow group. The program \textit{mori.x} only deals with simplicial toric varieties, 
for which the Picard group is always a finite index subgroup of the Chow group \cite{fulton}. Hence the Cartier divisors are integer
multiples of the Weil divisors and this ambiguity does not arise.

\subsection{The topological data of toric divisors} 

Using option ``-d'', \textit{mori.x} computes the arithmetic genera of the toric divisors restricted to the embedded hypersurface
and determines the del Pezzo candidates among them. The program checks the del Pezzo property against two necessary conditions:
First, for a del Pezzo divisor $S$ of type $n$, the following equations should hold:
\begin{equation}
\int_S c_1 (S)^2 = 9-n \, , \qquad \int_S c_2(S) = n+3 \qquad \Longrightarrow \qquad \chi_0(S)=\int_S \text{Td}(S)=1\, .
\end{equation} 
Here, \text{Td(S)} denotes the Todd class of $S$, which gives the zeroth arithmetic genus of $S$ upon integration. 
This test also allows to determine the type of the del Pezzo surface in question. 

The second necessary condition comes from the fact that a del Pezzo surface is a two--dimensional Fano manifold. 
Hence, the first Chern class of $S$ integrated over all curves on $S$ has to be positive:
\begin{equation}
 D_i \cap S\cap c_1(S) >0 \qquad\forall D_i:\, D_i\neq S \, ,\quad  \; D_i\cap S\neq 0 \, .
\end{equation}  
This condition would be sufficient if we were able to access \textit{all} curves of the hypersurface.
In our construction, however, we can only check for curves induced by toric divisors. This functionality was
added to carry out the analysis of base manifolds for elliptic fibrations in \cite{Knapp:2011wk}.

\section{Structure of the program}

\textit{mori.x} is part of the new releases of PALP (starting from version 2.0). 
The general structure of the package has not been changed, only some new files have been added.  
Hence, all general annotations to the package in \cite{Kreuzer:2002uu} remain valid.
In this section, we provide an overview of the composition of \textit{mori.x} and discuss its dependencies on pre--existing files.  

The source code of \textit{mori.x} is contained in the program files \textit{mori.c}, \textit{MoriCone.c}, \textit{SingularInput.c}, and the header file \textit{Mori.h}. 
\textit{Makefile} reflects the dependencies of the program. The compilation supports 32 as well as 64 bit architectures. One can adjust the compilation parameters according to ones needs. 
The optimization level is set at \verb+-03+ by default.   

\textit{mori.c} contains the main and the help information routines.
Further, basic manipulations (completion, calculation of facet equations,...) of the polyhedron are performed with the help of core routines from \textit{Vertex.c}.   
In particular, reflexivity of the polytope is checked.     

\textit{MoriCone.c} is at the heart of \textit{mori.x}. After determining the non--simplicial facets, 
their triangulation is performed by the routine \verb+GKZSubdivide+. 
This function identifies the maximal dimensional secondary fans of the facets and makes a case--by--case triangulation depending on their dimensions.
This function is only implemented for secondary fans up to dimension three.\footnote{The pre--compiler command TRACE\_TRIANGULATION in \textit{MoriCone.c} enables diagnostic information 
about the triangulation. This data might be of use for the motivated programmer who wants to extend the subdivision algorithm.}
Once the subdivision is accomplished, the program determines the Stanley--Reisner ideal (\verb+StanleyReisner+) and computes the Mori cone (\verb+Print_Mori+). 
Furthermore, it finds a basis of the toric divisor classes.

\textit{SingularInput.c} is the interface to SINGULAR \cite{DGPS}. 
The latter is a very efficient computer algebra system for computations with polynomial rings. 
In \textit{SingularInput.c}, the Chow ring is determined from the Stanley--Reisner ideal, the linear relations among the toric divisors, and a basis of the toric divisors of the triangulated polytope. 
This data is put together by \verb+HyperSurfSingular+ and redirected to SINGULAR, which then determines the intersection ring restricted to the hypersurface and computes its characteristic classes.
\footnote{The input for SINGULAR can be displayed in the standard output of \textit{mori.x} by turning on the pre--compiler definition TEST\_PRINT\_SINGULAR\_IO in \textit{SingularInput.c}.}

The most important routines of \textit{mori.x} are documented in the file \textit{Mori.h}. This header file provides a more detailed description of the structure of the program.



\begin{thebibliography}{9}

\bibitem{Kreuzer:2002uu}
  M.~Kreuzer, H.~Skarke,
  ``PALP: A Package for analyzing lattice polytopes with applications to toric geometry,''
  Comput.\ Phys.\ Commun.\  {\bf 157 } (2004)  87-106.
  [math/0204356 [math-sc]].

\bibitem{PALPwww}
	\url{http://hep.itp.tuwien.ac.at/~kreuzer/CY.html}


\bibitem{topcom}
 Jörg Rambau, ``TOPCOM: Triangulations of Point Configurations and Oriented Matroids'', 
 Mathematical Software - ICMS 2002 (Cohen, Arjeh M. and Gao, Xiao-Shan and Takayama, Nobuki, eds.), World Scientific (2002), pp. 330-340.

\bibitem{sage}
W.\thinspace{}A. Stein et~al., ``Sage Mathematics Software'', The Sage Development Team, \url{http://www.sagemath.org}


\bibitem{DGPS}
Decker, W.; Greuel, G.-M.; Pfister, G.; Sch{\"o}nemann, H.: 
``Singular {3-1-3} --- A computer algebra system for polynomial computations'' (2011).\\
\url{http://www.singular.uni-kl.de}



\bibitem{Odapark}
T.~Oda and H.~Park, ``Linear Gale transforms and Gelfand-Kapranov-Zelevinskij
decompositions,'' Thoku Math. J. {\bf 43} (1991) 375399.

\bibitem{Berglund:1995gd}
  P.~Berglund, S.~H.~Katz and A.~Klemm,
  ``Mirror symmetry and the moduli space for generic hypersurfaces in toric
  varieties,''
  Nucl.\ Phys.\  B {\bf 456} (1995) 153
  [arXiv:hep-th/9506091].


\bibitem{fulton}
W.~Fulton, ``Introduction to toric varieties'', Annals of mathematics studies, \textit{Princeton
University Press}, Princeton (1993).

T.~Oda,
  ``Convex Bodies and Algebraic Geometry,''
  {\it Springer-Verlag} (1988).


\bibitem{generaltoric}

  D.~A.~Cox, S.~Katz,
  ``Mirror symmetry and algebraic geometry,''
  Providence, USA: \textit{AMS} (2000) 469 p.
  
  M.~Kreuzer,
  ``Toric geometry and Calabi-Yau compactifications,''
  Ukr.\ J.\ Phys.\  {\bf 55 } (2010)  613.
  [hep-th/0612307].

  J.~Knapp, M.~Kreuzer,
  ``Toric Methods in F-theory Model Building,''
   [arXiv:1103.3358 [hep-th]].


\bibitem{Kreuzer:2000qv}
  M.~Kreuzer, H.~Skarke,
  ``Reflexive polyhedra, weights and toric Calabi-Yau fibrations,''
  Rev.\ Math.\ Phys.\  {\bf 14 } (2002)  343-374.
  [math/0001106 [math-ag]].


\bibitem{Batyrev:1994hm}
  V.~V.~Batyrev,
  ``Dual polyhedra and mirror symmetry for Calabi--Yau hypersurfaces in toric varieties,''
  J.\ Alg.\ Geom.\  {\bf 3 } (1994)  493-545.
  [alg-geom/9310003].



\bibitem{Billera} L. J. Billera, P. Filliman, B. Sturmfels, Constructions and complexity of secondary polytopes,
     Adv. in Math. {\bf 83} (1990) 155.

\bibitem{GKZ} I. M. Gelfand, M. M. Kapranov, A. V. Zelevinsky,Discriminants, Resultants, and Multidi-
mensional Determinants, \textit{Birkh\"auser}, Boston (1994).
                            

\bibitem{Knapp:2011wk}
  J.~Knapp, M.~Kreuzer, C.~Mayrhofer, N.-O.~Walliser,
  ``Toric Construction of Global F-Theory GUTs,''
  JHEP {\bf 1103 } (2011)  138.
  [arXiv:1101.4908 [hep-th]].

 


\end{thebibliography}
\end{document}